\documentclass[12pt,a4paper]{amsart}

\usepackage{latexsym, amssymb,amsmath, amsthm,amsfonts, mathrsfs}
\usepackage[arrow, matrix, curve]{xy}
\usepackage{hyperref}
\usepackage[short,nodayofweek]{datetime}
 \usepackage{verbatim}

\usepackage{nicefrac}

\newcommand{\ZZ}{\mathbb{Z}}
\newcommand{\NN}{\mathbb{N}}
\newcommand{\kk}{\Bbbk}
\newcommand{\kv}{{\kk[V]}}

\newcommand{\kvg}{{\kk[V]^G}}

\newcommand{\spec}{\mathrm{Spec}}
\newcommand{\done}{\hfill $\triangleleft$}

\newcommand{\isomto}{\overset{\sim}{\rightarrow}}

\newcommand{\GG}{\mathbb{G}}
\newcommand{\Ga}{{{\GG_a}}}
\newcommand{\kvga}{{\kv^{\Ga}}}
\newcommand{\kvd}{\kv^D}

\newcommand{\V}{{\mathcal{V}}}

\newcommand{\frakf}{\mathfrak{f}}

\def\SL2{\operatorname{SL}_{2}(K)}

\def\GL2{\operatorname{GL}_{2}(K)}
\def\Ga{{\mathbb G}_{a}}

\def\INVSL2{$K[V]^{operatorname{SL}_{2}(K)}$}
\def\INVSO2{$K[V]^{operatorname{SO}_{2}(K)}$}
\def\INVGL2{$K[V]^{operatorname{GL}_{2}(K)}$}

\def\GL{\operatorname{GL}}
\def\SL{\operatorname{SL}}

\newtheorem{lem}{Lemma}[section]
\newtheorem{thm}[lem]{Theorem}

\newtheorem{cor}[lem]{Corollary}

\theoremstyle{definition}

\theoremstyle{remark}
\newtheorem{rmk}[lem]{Remark}
\newtheorem{eg}[lem]{Example}

\newtheoremstyle{Acknowledgements}
  {}
    {}
     {}
     {}
    {\bfseries}
    {}
     {.5em}
     {\thmname{#1}\thmnumber{ }\thmnote{ (#3)}}
\theoremstyle{Acknowledgements}
\newtheorem{ack}{Acknowledgements.}

\title[Finite separating sets and quasi-affine quotients]{Finite separating sets and \\ quasi-affine  quotients}

\author{ Emilie Dufresne}
\address{Mathematisches Institut\\
Universit\"at Basel\\
Rheinsprung 21\\
4051 Basel, Switzerland}
\email{emilie.dufresne@unibas.ch}

\date{\today}



\begin{document}
\maketitle


\begin{abstract}
Nagata's famous counterexample to Hilbert's fourteenth problem shows that the ring of invariants of an algebraic group action on an affine algebraic variety is not always finitely generated. In some sense, however, invariant rings are not far from affine. Indeed, invariant rings are always quasi-affine, and there always exist finite separating sets. In this paper, we give a new method for finding a quasi-affine variety on which the ring of regular functions is equal to a given invariant ring, and we give a criterion to recognize separating algebras. The method and criterion are used on some known examples and in a new construction.
 \end{abstract}



\section{Introduction}

The ring of invariants of an algebraic group action on an affine variety is the subalgebra formed by those regular functions which are constant on the orbits. A central question in Invariant Theory, thought to be the inspiration for Hilbert's fourteenth problem, is to ask if the ring of invariants is always finitely generated, that is, if it is always equal to the ring of regular functions on some affine variety. Nagata  \cite{mn:ofph} gave a negative answer in 1958: a 32-dimensional linear representation of a non-reductive group. In 1990, Roberts   gave a new, significantly simpler counterexample: an action of the additive group on a 7-dimensional affine space (\cite{pr:aigsbpsrnchfp}, see Example \ref{eg-R7}). It lead to similar smaller examples by Freudenburg  in dimension 6 (\cite{gf:achfpd6}, see Example \ref{eg-F6}), and Daigle and Freudenburg  in dimension 5 (\cite{dd-gf:achfpd5}, see Example \ref{eg-DF5}), the smallest known counterexample to Hilbert's fourteenth problem.

Invariant rings are not far from finitely generated. Not only did Nagata   prove that they are at least rings of regular functions on some quasi-affine variety (see \cite[Chapter V.5]{mn:lfph}), but also Derksen and Kemper  showed that there always exists a finite separating set, that is, there always is a finite collection of invariants which can distinguish between any two points which are distinguished by some invariant (see \cite[Proposition 2.3.12]{hd-gk:cit}). The first result was made constructive by Derksen and Kemper in the case of an action of a connected unipotent group on a factorial variety (see  \cite[Algorithm 3.8]{hd-gk:ciagac}), but the algorithm is not very practical. The second result is highly non-constructive. Until now, only one example appeared in the literature: Winkelmann \cite{jw:irqq} found a quasi-affine variety on which the regular functions are the invariants from the Daigle-Freudenburg example (see Example \ref{eg-DF5}), and jointly with Kohls \cite{ed-mk:afssdfchfp}, we constructed a finite separating set for the same invariant ring.

In this paper, we give a new method for finding a quasi-affine variety on which the ring of regular functions is equal to a given invariant ring. In addition, we give a criterion to recognize separating algebras. The method and criterion are used on known examples in Section \ref{section-keg}, and to construct a new example in Section \ref{section-neg}.

\begin{ack}
Part of the work discussed here was done while I was a MATCH postdoctoral fellow in Heidelberg. I especially thank  Andreas Maurischat in Heidelberg, Gregor Kemper and Martin Kohls in Munich, Jonathan Elmer in Bristol, and Hanspeter Kraft in Basel.
\end{ack}


\section{Main Result}\label{section-main-result}

Let $\kk$ be an algebraically closed field, and let $G$ be an algebraic group over $\kk$. Suppose $G$ acts on $V$, an irreducible, normal affine algebraic $\kk$-variety (so that $\kv$, the ring of regular functions on $V$,  is a normal, finitely generated $\kk$-domain). Such an action induces a representation of $G$ on $\kv$ via $\sigma\cdot f=f\circ (-\sigma)$. The elements of $\kv$ which are fixed by $G$ form a subalgebra $\kvg$, called  the \emph{ring of invariants}.

By definition, invariants are constant on orbits. Thus, for two points $u,v\in V$ and $f\in\kvg$, if $f(u)\neq f(v)$, then $u$ and $v$ belong to distinct orbits, and we say $f$ \emph{separates} $u$ and $v$.  Accordingly, a subset $E\subseteq\kvg$ is called a \emph{separating set}  if any two points $u,v\in V$ which are separated by some invariant \cite[Definition 2.3.8]{hd-gk:cit} are separated by an element of $E$. A subalgebra $A\subset\kvg$ which is a separating set is called a \emph{separating algebra}. More generally, if $U$ is a subset of $V$, we say $E$ is a separating set on $U$ if the elements of $E$ separate any 2 points of $U$ which are separated by some invariant in $\kvg$ (cf. \cite[Definition 1.1]{gk:si}).

We recall the notation introduced in \cite[Section 2.1]{hd-gk:ciagac}, which fills the gap between colon operations on ideals and Nagata's ideal transform (see \cite[Chapter V.5]{mn:lfph}). If $A$ and $B$ are subsets of a commutative ring $S$, we define the following colon operations \cite[Definition 2.1]{hd-gk:ciagac}:
\[\begin{array}{l}
(A:B)_S:=\{f\in S\mid fB\subseteq A\}, \textrm{ and}\\
(A:B^\infty)_S:=\bigcup_{r=1}^\infty(A:B^r)_S=\{f\in S \mid \exists r\geq 1,\textrm{ such that }fB^r\subseteq A\}.\end{array}\]
When $A$ and $B$ are ideals, these are the usual colon ideals. If $R$ is a domain with field of fractions $Q(R)$ and $0\neq f\in R$, then 
\[(R:(f)^\infty)_{Q(R)}=R_f.\]
If $R=\kk[X]$ is the ring of regular functions on an irreducible affine variety and $Y$ is the zero set of the ideal $I$, then $(R:I^\infty)_{Q(R)}$ is the ring of regular functions on the quasi-affine variety $X\setminus Y$ \cite[Lemma 2.3]{hd-gk:ciagac}. Thus, $(R:I^\infty)_{Q(R)}$ corresponds exactly to the ideal transform of Nagata.

\begin{thm}\label{thm-trick}
Let $A\subset\kvg$ be a finitely generated subalgebra and let $f_1,\ldots,f_r\in A$ be such that $A_{f_i}=\kv^G_{f_i}$ for each $i$.
\begin{enumerate}
\item \label{thm-trick-1} If $\V_V(f_1,\ldots,f_r)\subseteq V$ has codimension at least 2 (that is, if $(f_1,\ldots,f_r)\kv$ has height at least $2$ in $\kv$), then $\kvg$ is equal to the ring of regular functions on the quasi-affine variety $\spec(A)\setminus \V(f_1,\ldots,f_r)$, that is, 
$ \kvg=(A:(f_1,\ldots,f_r)^\infty)_{Q(A)}$.
\item  \label{thm-trick-2} If $A$ is a  separating algebra on $\V_V(f_1,\ldots,f_r)\subseteq V$, then $A$ is a  separating algebra on all of $V$. 
\end{enumerate}

 \begin{proof}~\\
 (\ref{thm-trick-1}):~Our assumptions imply that $Q(A)=Q(\kvg)$. Since $\kv$ is normal and since $(f_1,\ldots,f_r)\kv$ has height at least 2, it follows that $(A:(f_1,\ldots,f_r)^\infty)_{Q(A)}$ is a subset of
\[ (\kv:(f_1,\ldots,f_r)^\infty)_{\kk(V)}\cap Q(A)=\kv\cap Q(\kvg)=\kvg.\]

Take $f\in \kvg$. Since $f\in\kv^G_{f_i}=A_{f_i}$, there is $s_i\geq 0$ such that $f_i^{s_i}f\in A$. If $s=s_1+\ldots +s_r$, then $f((f_1,\ldots,f_r)A)^s\subseteq A$. Therefore, $f\in(A:(f_1,\ldots,f_r)^\infty)_{Q(A)}$.
 
 \noindent
 (\ref{thm-trick-2}):~Suppose $u,v\in V$ are separated by $f\in\kvg$. If both $u$ and $v$ are in $\V_V(f_1,\ldots,f_r)$, our assumptions imply that $u$ and $v$ are separated by an element of $A$. If only one of $u,v$ is in $\V_V(f_1,\ldots,f_r)$, then $u$ and $v$ are separated by an $f_i$.  If neither $u$ nor $v$ is in $\V_V(f_1,\ldots,f_r)$, and if  no $f_i$ separates $u$ and $v$, then there is a $j$ such that $f_j(u)=f_j(v)\neq 0$. Since $\kvg=(A:(f_1,\ldots,f_r)^\infty)_{Q(A)}\subseteq A_{f_j}$, there exists $m\geq 0$ such that ${f_j}^mf\in A$. As $({f_j}^mf)(u)={f_j}(u)^mf(u)={f_j}(v)^mf(u)\neq {f_j}(u)^mf(u)=({f_j}^mf)(u)$, an element of $A$ separates $u$ and $v$.  \end{proof}
\end{thm}

If $B$ is a $\kk$-algebra, then
\[\frakf_B:=\{0\}\cup \{f\in B \mid B_f \textrm{ is a finitely generated }\kk \textrm{-algebra}\}\]
is a radical ideal of $B$, called the \emph{finite generation locus ideal} \cite[Proposition 2.10]{hd-gk:ciagac}. It is equal to $B$ exactly when $B$ is finitely generated. Using Theorem \ref{thm-trick} relies on finding enough elements in the finite generation ideal.

\begin{cor}\label{cor-sep}
Suppose $A$ and $f_1,\ldots,f_r$ satisfy the conditions of Theorem \ref{thm-trick}(\ref{thm-trick-1}). If $\kvg\subseteq \kk + (f_1,\ldots,f_r)\kv$, then $A$ is a separating algebra.
\begin{proof}
If if $\kvg\subseteq \kk + (f_1,\ldots,f_r)\kv$, then all invariants are constant on $\V_V(f_1,\ldots,f_r)$ and so the condition of Theorem \ref{thm-trick}(\ref{thm-trick-2}) is automatically satisfied.
\end{proof}
\end{cor}

\begin{rmk}
Theorem \ref{thm-trick} implies that if we have $A\subset\kvg$ and an ideal $I$ of $A$ such that $\kvg$ is equal to the ring of regular functions on $\spec(A)\setminus \V(I)$, then under some additional assumptions, $A$ is a separating algebra. On the other hand, if $A\subseteq\kvg$ is a normal finitely generated separating algebra with $Q(A)=Q(\kvg)$, then there is an ideal $I$ of $A$ such that $\kvg$ is equal to the ring of regular functions on $\spec(A)\setminus \V(I)$. This can be deduced from \cite[Theorem 2 and Lemma 7]{jw:irqq} as follows.

Each  $E\subseteq\kv$ induces an equivalence relation $\sim_E$ on $V$. For $u,v\in V$, we write $u\sim_Ev$ if and only if $f(u)=f(v)$ for all $f\in E$.
In \cite[Lemma 7]{jw:irqq}, Winkelmann shows that there exists a normal, finitely generated subalgebra $A\subset \kvg$ such that $\sim_A=\sim_{\kvg}$. In the proof of \cite[Theorem 2]{jw:irqq}, he shows that we can assume $A$ is normal and $Q(A)=Q(\kvg)$. Then, there is  an ideal $I$ of $A$ such that
\[\kvg=(A:I^\infty)_{Q(A)}.\]
The key observation is that $\sim_A=\sim_{\kvg}$ exactly when $A$ is a  separating algebra. In particular, \cite[Lemma 7]{jw:irqq} implies the existence of a finitely generated separating algebra.
\end{rmk}


\section{Additive group actions}\label{section-additive-group-actions}

 For the examples discussed in Sections \ref{section-keg} and \ref{section-neg}, we concentrate on algebraic actions of the additive group $\Ga=(\kk,+)$, and assume that $\kk$ has characteristic 0. Such an action corresponds to a \emph{locally nilpotent derivation (LND)}, that is, a $\kk$-linear map $D:\kv\rightarrow \kv$ such that 
\begin{enumerate}
\item for all $a,b\in \kv$, we have $D(ab)=aD(b)+bD(a)$, and
\item for all $b\in \kv$, there exists $m\geq 0$ such that $D^m(b)=0$.
\end{enumerate}
The $\Ga$-action on $V$ is given by the $\kk$-algebra homomorphism:
\[\begin{array}{lccc}
\theta : & \kv & \longrightarrow & \kv\otimes_\kk\kk[T]\\
            & f     & \longmapsto     & \theta(f),\end{array}\]
 where $\kk[T]=\kk[\Ga]$ is the ring of regular functions on the algebraic group $\Ga$. This $\Ga$-action induces an action on $\kv$ via $a\cdot f=\theta(f)|_{T=-a}$. The correspondence between $D$ and the $\Ga$-action is given by 
\[\theta(f)=\sum_{k=0}^\infty\frac{D^k(f)}{k!}T^k.\]
The ring of invariants $\kv^{\Ga}$ coincides with the kernel of $D$, which we write $\kv^D$. For convenience, we will describe $\Ga$-actions on $V$ by giving the corresponding LND on $\kv$. For more information on LND, we refer to the excellent book of Freudenburg \cite{gf:atlnd}. 

If $s\in \kv$ is a slice, that is, if $D(s)=1$, then $s:V\rightarrow\Ga$ is $\Ga$-equivariant, and there is a $\Ga$-equivariant isomorphism $\Ga\times s^{-1}(0){\isomto}V$, given by $(a,v)\mapsto a\cdot v$, which identifies the invariants $\kv^{\Ga}$ with $\kk[s^{-1}(0)]$. In general, if $s\in\kv$, $f:=D(s)\neq 0$, and $D^2(s)=0$, then $s/f$ is a slice on $V_f=V\setminus \V_V(f)$, where $\V_V(f)$ denotes the zero set of $f$ in $V$. Such an $s$ is called a \emph{local slice}. We then obtain generators for $\kv^D_f$ as follows (it is the first step of van den Essen's algorithm):

\begin{lem}[see {\cite[Sections 3 and 4]{ae:aacigaav}}]\label{lem-slice}
Take $s\in \kv$ such that $f=D(s)\neq 0$ and $D^2(s)=0$. If $\kv=\kk[b_1,\ldots,b_r]$, then $\kv^D_f$ is generated by $f$, $1/f$, and $\{f^{e_i}\theta(b_i)|_{T=-{s}/{f}}\mid i=1,\ldots,r\}$, where $e_i$ is minimal so that $f^{e_i}\theta(b_i)|_{T=-{s}/{f}}\in\kv$.
\end{lem}

As $\Ga$ is a connected unipotent group, when $V$ is factorial, the finite generation ideal $\frakf_{\kvd}$ generates an ideal of height at least 2 in $\kv$ (see the proof of correctness of \cite[Algorithm 2.22]{hd-gk:ciagac}). Therefore, there always exist $\{f_1,\ldots,f_m\}$ and $A$ satisfying the conditions of Theorem \ref{thm-trick}(\ref{thm-trick-1}). Combining some existing algorithms, one can compute such $\{f_1,\ldots,f_m\}$ and $A$ as follows. First, use \cite[Algorithm 3.20]{tk:acir} to compute $f_0\in\frakf_{\kvd}$, and $g_{0,1},\ldots,g_{0,s_0}\in\kv^N$, such that $\kvd_{f_0}=\kk[f_0,1/f_0,g_{0,1},\ldots,g_{0,s_0}]$. Next, use \cite[Algorithm 2.13]{hd-gk:ciagac} with $S=\kv$, $R_0=\kk[f_0,g_{0,1},\ldots,g_{0,s_0}]$, and $\mathfrak{a}=f_0R$ to compute further elements $\{f_1,\ldots,f_r\}$  of $\frakf_{\kvd}$ until the ideal $(f_0,f_1,\ldots,f_r)\kv$ has height at least 2. 
The last step is to use van den Essen's Algorithm to compute $g_{i,1},\ldots,g_{i,s_i}\in\kvd$ such that $\kvd_{f_i}=\kk[f_i,1/f_i,g_{i,1},\ldots,g_{i,s_i}]$. Taking $\{f_0,\ldots,f_r\}$ and $A=\kk[f_0,\ldots,f_r,g_{i,j}\mid i=0,\ldots,r, j=1,\ldots,s_r]$ will satisfy the conditions of Theorem \ref{thm-trick}(\ref{thm-trick-1}).


\section{First examples}\label{section-keg}


\begin{eg}[Daigle and Freudenburg \cite{dd-gf:achfpd5}] \label{eg-DF5}
 Let $V:=\kk^5$, and let  $R:=\kk[x,s,t,u,v]$ be the ring of regular functions on $V$. Define a LND on $R$ via:
 \[\Delta:=x^3\frac{\partial}{\partial s}+s\frac{\partial}{\partial t}+t\frac{\partial}{\partial u}+x^2\frac{\partial}{\partial v}.\]
 
 Daigle and Freudenburg proved in  \cite{dd-gf:achfpd5} that the ring of invariants $R^{\Delta}$ is not finitely generated. In \cite[Section 4]{jw:irqq}, Winkelmann defined a subalgebra
 \[\begin{array}{rl}     
 A:=&\kk[f_1,f_2,f_3,f_4,f_5,f_6]\footnotemark\\  
    =&\kk[x, 2x^3t-s^2, 3x^6u-3x^3ts+s^3, xv-s, x^2ts-s^2v+\\
 &\hspace{.7cm} 2x^3tv-3x^5u, -18x^3tsu+9x^6u^2+8x^3t^3+6s^3u-3t^2s^2],
 \end{array}\]  
  and proved that $R^\Delta$ is equal to the ring of regular functions on $\spec(A)\setminus \V(x,2x^3t-s^2)$.
 With Kohls \cite{ed-mk:afssdfchfp}, we proved that $A$  is a separating algebra. We will show how both results follow from Theorem \ref{thm-trick}.
    \footnotetext{ In fact, we have $A=\kk[f_1,f_2,f_4,f_5,f_6]$, since $f_3=-f_1f_5+f_2f_4$.}
 
   We have $x^3=\Delta(s)\in R^\Delta$ and $f_2=2x^3t-s^2=\Delta(3x^3u-st)\in R^\Delta$. Lemma \ref{lem-slice} yields the following generators 
   \[
   \begin{array}{ll}
   \textrm{for }R_x^\Delta=R_{x^3}^\Delta:                                                 &    \textrm{and for }R_{f_2}^\Delta:\\
   x, 1/x v                                                                                                      &    f_2,1/{f_2}\\
   \theta(x)|_{{-s}/{x^3}}=x=f_1                                                        & \theta(x)|_{{-(3x^3u-st)}/{f_2}}=x=f_1\\
   \theta(s)|_{{-s}/{x^3}}=0                                                                &   f_2\theta(s)|_{{-(3x^3u-st)}/{f_2}}=-f_3/2\\
   x^3\theta(t)|_{{-s}/{x^3}}=f_2/2                                             &  f_2^2\theta(t)|_{{-(3x^3u-st)}/{f_2}}=x^3f_6/2\\
   x^6\theta(u)|_{{-s}/{x^3}}=f_3/3\ \hspace{2cm} & f_2^3\theta(u)|_{{-(3x^3u-st)}/{f_2}}=f_6f_3/6\\
   x\theta(v)|_{{-s}/{x^3}}=f_4,                                                            & f_2\theta(v)|_{{-(3x^3u-st)}/{f_2}}=f_5.\\
   \end{array}\]
  Observe that $A$ contains the above polynomials, and so $A_x=R_x^\Delta$ and $A_{f_2}=R_{f_2}^\Delta$. As $\V_{\kk^5}(x^3,2x^3t-s^2)=\V_{\kk^5}(x,s)$ has codimension 2, Theorem \ref{thm-trick}(\ref{thm-trick-1}) implies that $R^\Delta$ is the ring of regular functions on $\spec(A)\setminus\V(x,2x^3t-s^2)$. 

Using the fact that $\Delta$ is graded and commutes with $\frac{\partial}{\partial u}$ and  $\frac{\partial}{\partial v}$, one can show that ${R}^{\GG_a}\subseteq \kk\oplus(x,s){R}$ (see \cite[Proposition 3.2]{ed-mk:afssdfchfp}).  Corollary \ref{cor-sep} then imply that $A$ is a separating algebra.\done\\
  \end{eg}


\begin{eg}[Roberts \cite{pr:aigsbpsrnchfp}] \label{eg-R7}
Let $B:=\kk[x_1,x_2,x_3,y_1.y_2,y_3,v]$, and let  $2\leq m\in\ZZ$. Consider the LND defined on $B$ via:
\[D:=x_1^{m+1}\frac{\partial}{\partial y_1}+ x_2^{m+1}\frac{\partial}{\partial y_2}+ x_3^{m+1}\frac{\partial}{\partial y_3}+ (x_1x_2x_3)^{m}\frac{\partial}{\partial v}.\]
Roberts \cite{pr:aigsbpsrnchfp} proved that $B^D$ is not finitely generated.  

For each $i$, $D(y_i)=x_i^{m+1}\in B^D$, and Lemma \ref{lem-slice} yields the following invariants:

{\footnotesize \[\begin{array}{lll}
\phi_1  = x_1^{m+1}y_2-x_2^{m+1}y_1, &
\phi_2  = x_1^{m+1}y_3-x_3^{m+1}y_1, &
\phi_3  = x_2^{m+1}y_3-x_3^{m+1}y_2, \\
\phi_4  =  (x_1x_2)^my_3-x_3v, &
\phi_5  =  (x_1x_3)^my_2-x_2v, &
\phi_6  =  (x_2x_3)^my_1-x_1v.
  \end{array}\]}Let $A=\kk[x_1,x_2,x_3,\phi_1,\phi_2,\phi_3,\phi_4,\phi_5,\phi_6]$. By construction, we have  $B^D_{x_i}=A_{x_i}$. As $\V_{\kk^7}(x_1,x_2,x_3)$ has codimension 3, Theorem \ref{thm-trick}(\ref{thm-trick-1}) implies that $B^D$ is the ring of regular functions on $\spec(A)\setminus\V(x_1,x_2,x_3)$.

As $B^{D}\subseteq \kk\oplus(x_1,x_2,x_3)B$ (see \cite[Lemma 2]{pr:aigsbpsrnchfp}),   Corollary  \ref{cor-sep} imply $A$ is a  separating algebra.\done\\
\end{eg}

\begin{rmk}
Part (\ref{thm-trick-1})  of Theorem \ref{thm-trick} does not imply part (\ref{thm-trick-2}). Indeed,  $A':=\kk[x_1, x_2, x_3, \phi_1, \phi_2, \phi_3, \phi_5, \phi_6]$ is not a separating algebra, as $A'$ does not separate $(0,0,1,0,0,0,1)$ from the origin, although $\phi_4(0,0,1,0,0,0,1)=1$. On the other hand, by Lemma \ref{lem-slice},  $A'_{x_1}=B^D_{x_1}$ and $A'_{x_2}=B^D_{x_2}$. Since $\V_{\kk^7}(x_1,x_2)$ has codimension 2, Theorem \ref{thm-trick}(\ref{thm-trick-1}) implies that $B^D$ is the ring of regular functions on $\spec(A')\setminus \V(x_1,x_2)$.
\end{rmk}


\begin{eg}[Freudenburg \cite{gf:achfpd6}] \label{eg-F6}
Let $B:=\kk[x,y,s,t,u,v]$, and define a LND on $B$ via:
\[D:=x^3\frac{\partial}{\partial s}+y^3s\frac{\partial}{\partial t}+y^3t\frac{\partial}{\partial u}+x^2y^2\frac{\partial}{\partial v}.\]
 Freudenburg \cite{gf:achfpd6} showed that $B^D$ is not finitely generated.

 Let $A$ be the $\kk$-algebra generated by:
\[\begin{array}{l}
x,~    y,~    -y^2s + xv,~    -\frac{1}{2}y^3s^2 + x^3t,\\
    -x^2y^3st + 3x^5u + y^4s^2v - 2x^3ytv,\\
    -\frac{3}{2}y^6s^2t^2 + 4x^3y^3t^3 + 3y^6s^3u - 9x^3y^3stu + \frac{9}{2}x^6u^2.
\end{array}\]
 We have $D(s)=x^3\in B^D$ and $D(3x^3u-y^3st)=2x^3y^3t-y^6s^2\in B^D$. Comparing with with the generators given by  Lemma \ref{lem-slice}, we see that $A_{2x^3y^3t-y^6s^2}=B^D_{2x^3y^3t-y^6s^2}$, and $A_x=B^D_x$. As 
 \[\V_{\kk^6}(x,2x^3y^3t-y^6s^2)=\V_{\kk^6}(x,ys)=\V_{\kk^6}(x,y) \cup \V_{\kk^6}(x,s)\]
  has codimension 2, Theorem  \ref{thm-trick}(\ref{thm-trick-1}) implies that $B^D$ is the ring of regular functions on $\spec(A)\setminus\V(x,2x^3y^3t-y^6s^2)$

  We have $B^D\subseteq\kk\oplus(x,y)B$ (see \cite[Lemma 1]{gf:achfpd6}). A careful study of the list of generators given by Tanimoto  for $B^D$ \cite[Theorem 1.6]{rt:ofcfph} reveals that $B^D\subseteq \kk[y]\oplus(x,s)B$. As $A$ contains $y$, $A$ is a separating algebra on both $\V_{\kk^6}(x,y)$ and $\V_{\kk^6}(x,s)$. Hence, it is a separating algebra on $\V_{\kk^6}(x,2x^3y^3t-y^6s^2)$. By Theorem \ref{thm-trick}(\ref{thm-trick-2}), $A$ is a separating algebra on  all of $\kk^6$.  
\end{eg}


\section{A new 7-dimensional example}\label{section-neg}

The new 7-dimensional example discussed in this section illustrates the difficulty involved in applying Theorem \ref{thm-trick}(\ref{thm-trick-2}).

 Let $B:= \kk [x_1,x_2,x_3,y_1,y_2,y_3,v]$, and define a LND on $B$ via:
 \[D:=x_1^{a}\frac{\partial}{\partial y_1}+ x_2^{a}\frac{\partial}{\partial y_2}+ x_3^{a}\frac{\partial}{\partial y_3}+ (y_1y_2y_3)^{b}\frac{\partial}{\partial v},\]
  where $1\leq a,b\in\ZZ$.  We do not know if $B^D$ is finitely generated.

Noting that $D(y_i)=x_i^a\in B^D$, we apply  Lemma \ref{lem-slice} to define $A\subset B^D$ so that $A_{x_i}=B^D_{x_i}$:
 $$A:=\kk[x_1,x_2,x_3,x_1^ay_2-x_2^ay_1,x_1^ay_3-x_3^ay_1,x_2^3y_3-x_3^ay_2,h_1,h_2,h_3],$$ 
 where
\[h_i=x_i^{(2b+1)a}\theta(v)|_{T={-y_i}/{x_i^a}},~i=1,2,3,\] 
and $\theta:B\rightarrow B[T]$ is the map giving the $\Ga$-action. As $\V_{\kk^7}(x_1,x_2,x_3)$ has codimension 3, by Theorem \ref{thm-trick}(\ref{thm-trick-1}) $B^D$ is the ring of regular functions on $\spec(A)\setminus \V(x_1,x_2,x_3)$.

In Lemma \ref{lem-finalstep} below, we will show that $B^D\subseteq\kk\oplus (x_1,x_2,x_3)B$. By Corollary \ref{cor-sep}, it then follows that $A$ is a separating algebra.


Our argument to prove Lemma \ref{lem-finalstep} relies on the relationship between our new  7-di\-men\-sional example and a generalization of an example first proposed by Maubach \cite[Chapter 5]{sm:pekd}. Let $R:=\kk[x,y,z,u,w]$, and let $1\leq a,b \in \ZZ$. Define a LND on $R$:
 \[\Delta:=x^a\frac{\partial}{\partial y}+y\frac{\partial}{\partial z}+z\frac{\partial}{\partial u}+u^b\frac{\partial}{\partial w}.\]
In the case $a=1,b=2$, Maubach asked if $R^\Delta$ is finitely generated. The question remains open. 

In \cite[Section 7.2.3]{gf:atlnd}, Freudenburg explains how the Dai\-gle-Freu\-den\-burg example (see Example \ref{eg-DF5}) can be derived from Roberts's example (see Example \ref{eg-R7}) by ``removing all symmetries''. We follow the same argument to derive Maubach's example from our new 7-dimensional example.  Consider the faithful action on $B$ by the $3$-dimensional multiplicative group $\GG_m^3$ given by:
\[\begin{array}{l}
(\lambda,\mu,\nu)\cdot (x_1,x_2,x_3,y_1,y_2,y_3,v):=\\
\hspace{4cm}(\lambda x_1,\mu x_2, \nu x_3, \lambda^a y_1, \mu^a y_2, \nu^a y_3, (\lambda\mu\nu)^bv).
\end{array}\]
This action commutes with $D$. Additionally, $D$ commutes with the action of the symmetric group  $S_3$ given by:
\[\sigma\cdot (x_1,x_2,x_3,y_1,y_2,y_3,v):=(x_{\sigma(1)},x_{\sigma(2)},x_{\sigma(3)},y_{\sigma(1)},y_{\sigma(2)},y_{\sigma(3)},v).\]
The group $S_3$ acts  on $\GG_m^3$ by conjugation, and so $\GG_m^3\rtimes S_3$ acts on $B$. Since the $\GG_m^3$-action has no non-constant  invariants, we consider the subgroup $H$ of $\GG_m^3$ given by $\lambda\mu\nu=1$. This subgroup $H$ is isomorphic to $\GG_m^2$, and the group $G:=H\rtimes S_3$ acts on $B$. The invariant ring of $H$ is generated by monomials:
\[B^H=\kk[x_1x_2x_3, x_1^ay_2y_3, x_2^ay_1y_3, x_3^ay_1y_2, x_1^ax_2^ay_3, x_1^ax_3^ay_2, x_2^ax_3^ay_1, y_1y_2y_3, v].\]
Since $H$ is normal in $G$, $B^G=(B^H)^{S_3}$. Moreover, $B^G$ is a polynomial ring in 5 variables given as a subalgebra of $B$ by:
\[\kk[x_1x_2x_3, x_1^ay_2y_3+x_2^ay_1y_3+x_3^ay_1y_2, x_1^ax_2^ay_3+x_1^ax_3^ay_2+x_2^ax_3^ay_1, y_1y_2y_3, v].\]
Setting
\[\begin{array}{l}
   x:=x_1x_2x_3,\\
   y:=(x_1^ax_2^ay_3+x_1^ax_3^ay_2+x_2^ax_3^ay_1)/3,\\
   z:=(x_1^ay_2y_3+x_2^ay_1y_3+x_3^ay_1y_2)/6,\\
   u:=y_1y_2y_3/6,\\
   w:=v/6^b,
  \end{array}\]
  we have $B^G=R$, and the LND induced by $D$ coincides with $\Delta$. As $G$ is a reductive group and since $R^\Delta=(B^D)^G$, if $B^D$ is finitely generated, so is $R^\Delta$.


\begin{lem}\label{lem-step1}
$B^D \subseteq \kk[y_1,y_2,y_3]\oplus(x_1,x_2,x_3)B$.
\begin{proof}
If $B':=\kk[y_1,y_2,y_3,v]\cong B/(x_1,x_2,x_3)$, then $D$ induces a locally nilpotent derivation on $B'$:
\[D':=(y_1y_2y_3)^b\frac{\partial}{\partial v},\]
with kernel ${B'}^{D'}=\kk[y_1,y_2,y_3]$. Thus, if $f\in\kvga$, we can write
\[f:=x_1f_1+x_2f_2+x_3f_3+h,\]
where $h$ can be viewed as an element of $B'$. As $D(f)=0$, we have $D'(h)=0$, and so $f\in \kk[y_1,y_2,y_3]\oplus(x_1,x_2,x_3)B$. \end{proof}
\end{lem}

\begin{lem}\label{lem-step2}
$R^\Delta\subseteq \kk\oplus(x,y,z)R$.
\begin{proof}
If $R':=\kk[y,z,u,w]\cong R/(x)$, then $\Delta$ induces a locally nilpotent derivation on $R'$:
\[\Delta':=y\frac{\partial}{\partial z}+z\frac{\partial}{\partial u}+u^b\frac{\partial}{\partial v}.\]
Since this is an elementary monomial derivation in four variables, the ring of invariants is generated by at most four elements \cite{sm:tmdkkgm4e}, which we compute with van den Essen's Algorithm \cite[Section 4]{ae:aacigaav}. First, we write down the algebra map $\theta':R'\rightarrow R'[T]$ corresponding to $\Delta'$:
\[\begin{array}{l}
\theta'(y)=y,\\
   \theta'(z)=z+yT,\\
   \theta'(u)=u+zT+\frac{1}{2}yT^2,\\
\theta'(w)=\\ \hspace{0.5cm}w+\sum_{l=1}^{2b+1}\frac{1}{l}T^l
         \sum_{m=0}^{b}\binom{b}{b-m,2m+1-l,l-m-1}
                 \frac{1}{2^{l-m-1}}u^{b-m}z^{2m+1-l}y^{l-m-1}
\end{array}\]
Choosing the local slice $z$, the first step of the algorithm yields the following three generators:
\[ \begin{array}{l}
   y,\\
   h:=yu-\frac{1}{2}z^2,\\
   h':=\\ \hspace{0.3cm}y^{b+1}w+\sum_{l=1}^{2b+1}\frac{(-1)^l}{l}
         \sum_{m=0}^{b}\binom{b}{b-m,2m+1-l,l-m-1}
                \frac{1}{2^{l-m-1}}u^{b-m}z^{2m+1}y^{b-m}. 
  \end{array}\]
The second step of the algorithm yields the fourth generator:
\[h''=\frac{1}{y^n}\left(\frac{1}{\alpha^2}h^{2b+1}+2^{2b+1}h'^{2}\right),\]
where 
\[\alpha:=\sum_{j=0}^{b}\frac{(-1)^{j+b+1}}{(j+b+1)2^{j}}\binom{b}{j}=\nicefrac {(-1)^{b+1}2^b}{\binom{2b+1}{b+1}},\]
 and $n$ is maximal so that $y^n$ divides $\frac{4}{3}h^2+\frac{1}{\alpha}h'^{2b+1}$. It only remains to check that $y,h,h',h''\in(y,z)R'$. This is clear for $y$, $h$, and $h'$. Modulo $z$, we have
\[h''\equiv\frac{1}{y^n}\left(\frac{1}{\alpha^2}(uy)^{2b+1}+2^{2b+1}(wy^{b+1})^2\right).\]
Since the terms divisible by $y^{2b}$ in $\frac{1}{\alpha^2}{h}^{2b+1}$ and $2^{2b+1}{h'}^{2}$ do not cancel, $n\leq2b$. It follows that $h''\in(y,z)R'$, and so $R'^{\Delta'}\subseteq\kk\oplus(y,z)R'$, hence $R^\Delta\subseteq \kk\oplus(x,y,z)R$.
\end{proof}\end{lem}

\begin{lem}\label{lem-finalstep}
$B^D\subseteq \kk \oplus (x_1,x_2,x_3)B$.
\begin{proof}

The linear $\GG_m^3$-action on $B$ induces a $\ZZ^3$-grading (here, really a $\NN^3$-grading) $\omega$ on $B$ via characters where $B_{(0,0)}=B^G=\kk$ (see \cite[Proposition 4.14]{sm:aiim}). The derivation $D$ commutes with the $\GG_m^3$-action, implying that $B^D\subset B$ is a $\NN^3$-graded subalgebra. Therefore, it suffices to show that  every non-constant $\omega$-homogeneous element of $B^D$ is in the ideal $(x_1,x_2,x_3)B$.

Suppose, for a contradiction, that $f$ is a non-constant $\omega$-homogeneous element of $B^D$ not contained in the ideal $(x_1,x_2,x_3)B$. By Lemma \ref{lem-step1}, it is of the form $f=f_1+f_2$, where $f_1\in\kk[y_1,y_2,y_3]$ and $f_2\in(x_1,x_2,x_3)B$. As $f$ is $\omega$-homogeneous, so are $f_1$ and $f_2$. Hence, $f_1$ is supported at the monomial $y_1^{l_1}y_2^{l_2}y_3^{l_3}$. Let $F$ be the orbit product of $f$ under the $S_3$-action. We then have $F=F_1+F_2$, where $F_2\in(x_1,x_2,x_3)B$ and $F_1\in\kk[y_1,y_2,y_3]$ is supported at the monomial $y_1^ly_2^ly_3^l$. As $D$ commutes with the $S_3$-action, $D(F)=0$. The linear action of $H\cong \GG_m^2$ induces a $\ZZ^2$-grading on $B$ via characters, where $B_{(0,0)}=B^H$ (see again \cite[Proposition 4.14]{sm:aiim}). Let $F'$ be the component of $F$ of degree $(0,0)$, then $F'$ is $H$-invariant and contains the term $F_1$. As $S_3$ acts on $B^H$ and $F$ is $S_3$-invariant, $F'$ is $S_3$-invariant. It follows that $F'\in (B^H)^{S_3}=B^G=R$.  As $D$ commutes with the $H$-action, $D$ is graded with respect to the induced $\ZZ^2$-grading, and so $D(F')=0$, that is, $F'$ is an element of $R^\Delta$ containing supported at the monomial $u^l$, a contradiction to Lemma~\ref{lem-step2}.
\end{proof}\end{lem}

\begin{rmk}
As in our joint work with Maurischat \cite{ed-am:otgagipc}, one can define a characteristic-free analog to this new 7-dimensional example. The map $\theta$ has rational coefficients with denominators all dividing $(3b+1)!$. Thus, we can interpret $\theta$ as a locally finite iterative higher derivation over any field of characteristic $p>3b+1$. Use \cite[Theorem 1.1]{rt:aacklfihd} (the positive characteristic analog of Lemma \ref{lem-slice}) to define $A\subset B^\theta$ so that $A_{x_i}=B^D_{x_i}$:
\[A=\kk[x_1,x_2,x_3,x_1^ay_2-x_2^ay_1,x_1^ay_3-x_3^ay_1,x_2^ay_3-x_3^ay_2,f_1,f_2,f_3],\]
 where
\[f_i=x_i^{(2b+1)a}\theta(v)|_{T=\frac{-y_i}{x_i^a}},~i=1,2,3.\] 
 Theorem \ref{thm-trick}(\ref{thm-trick-1}) implies that $B^D$ is the ring of regular functions on $\spec(A)\setminus\V(x_1,x_2,x_3)$. We can  show that $B^\theta\subseteq \kk\oplus (x_1,x_2,x_3)B$, and so, by Corollary \ref{cor-sep}, $A$ is a separating algebra. The only significant difference with the characteristic zero case is that in Lemma \ref{lem-step2}, we must prove that the algorithm really ends after obtaining the fourth generator. This can be done as in the original argument of Maubach \cite[Case 3, theorem 3.1]{sm:tmdkkgm4e}, using that modulo $y$, $h''$ does not depend only on $z$.
\end{rmk}


\bibliographystyle{plain}
\bibliography{reference}
\end{document}